\documentclass[12pt, twoside]{article}
\usepackage{amsmath,amssymb}
\usepackage{amsthm}
\usepackage{graphicx}
\usepackage{epstopdf}
\usepackage{epsfig}
\usepackage{cite}
\usepackage{subfigure}
\usepackage[colorlinks,linkcolor=blue,anchorcolor=blue,citecolor=blue]{hyperref}
 \usepackage{color}

\newcommand{\xiaosanhao}{\fontsize{14pt}{\baselineskip}\selectfont}

\begin{document}

\title{\textbf{The generalized Bernoulli numbers and its relation with the Riemann zeta function at odd-integer arguments }}
\author{
Yayun Wu{\small\textsc{\thanks{E-mail address:~wyyfde@sina.com (Y. Wu).} }}\\
School of Mathematics and Statistics, Hefei Normal University,\\
Hefei, 230601, People's Republic of China}

\date{}
\maketitle
\textbf{Abstract:}
By using the generalized Bernoulli numbers, we deduce new integral representations for the Riemann zeta function at positive odd-integer arguments. The explicit expressions enable us to obtain criteria for the dimension of the vector space spanned over the rational by the $\zeta(2n+1)/\pi^{2n}$, $n\geq1$.



\section*{\xiaosanhao 1. Introduction }
The Riemann zeta function $\zeta(s)$ is defined for complex values of $s$, $\Re(s)>1$  as $\zeta(s):=\sum_{v=1}^{\infty}v^{-s}$, and its meromorphic continuation over the whole complex $s-$plane, except for a simple pole at $s=1$.  The research of the arithmetical properties of the Riemann zeta function at positive integers is one of the fascinating topics of the complex analysis and number theory. Finding recurrence formulas and integral representations of zeta function has become a significant issue. As is well known, for Riemann's zeta function $\zeta(s)$ at even points, we have the formula (found by Euler)
\begin{align}
\zeta(2k)=\frac{(-1)^{k-1}(2k)^{2k}B_{2k}}{2(2k)!},~ k\in \mathbb{N}^{+}.\nonumber
\end{align}
Here the $B_{k}=B_{k}(0)$ is the well-tabulated Bernoulli numbers defined by means of the generating function
\begin{align}
\frac{ze^{zx}}{e^{z}-1}=\sum_{k=0}^{\infty}\frac{B_{k}(x)}{n!}z^{k},~|z|<2\pi.\nonumber
\end{align}

However, there is no analogous closed evaluation for the values of $\zeta(s)$ at positive odd integers or fractional points up to the present time. Using the series of special functions and the multidimensional  integral to investigate $\zeta(s)$ are the main methods.
Connecting Euler polynomials and trigonometric functions, Srivastava and Choi[3] give the Riemann zeta function at positive odd integers
\begin{align}
\zeta(2n+1)=\frac{(-1)^{n}\pi^{2n+1}}{4\delta(1-2^{-(2n+1)})(2n)!}\int_{0}^{\delta}E_{2n}(u)\csc(u\pi){\rm d}u,~n\in \mathbb{N}^{+},~\delta=1~{\rm or}~1/2,\nonumber
\end{align}
where $E_{n}(u)$ is the Euler polynomials. Zurab[1], Cvijovi$\acute{c}$ and Klinowski[2] establish similar results by using Euler polynomials and the polylogarithm function[8]. Generalizing the trigonometric function to hyperbolic version, Lima[6], D'Avanzo and Krylov[5] research the integral result about $\zeta(n)$,
\begin{align}
\zeta(n)=\frac{2^{n}}{(2^{n}-1)(n-1)!}\int_{0}^{\infty}\ln^{n-1}(\coth(x)){\rm d}x,~n\geq2.\nonumber
\end{align}

These integral representations about $\zeta(s)$ hardly show the correlation with the irrational problem of the zeta function values at odd integers. Since the problem of the irrationality of the values of the zeta function at odd integers is one of the most attractive topics in modern number theory.
These problems motivate us to modify the integral representations about $\zeta(s)$.

The article is organized as follows. In $\S$ 2, we reconsider the evaluation of $\zeta(3)$ so that technical details of the general case do not obscure. In $\S$ 3, we elaborate the general case by using the generalized Bernoulli polynomials and numbers.
In $\S$ 4, we evaluate the numerical results about $\zeta(5)$ and $\zeta(7)$. In $\S$5, the formula for $\zeta(2n+1)$ enables us to obtain  linear forms with rational coefficients in the values of $\zeta(2n+1)/\pi^{2n}$, $n\in \mathbb{N}^{+}$.
We also ponder over the arithmetical properties of the sequence $\zeta(3)/\pi^{2},\zeta(5)/\pi^{4},\zeta(7)/\pi^{6}$,$\cdots$.

\section*{2. Evaluation of $\zeta(3)$}

In this part, we rediscover an integral representation for $\zeta(3)$. Let us briefly describe the proof for a better exposition of the general case and its consequences.

{\bf Lemma~2.1}
{\it For Riemann zeta function value $\zeta(3)$, we have \begin{align}
\zeta(3)=\frac{2\pi^{2}}{7}\int_{-\infty}^{+\infty}\frac{e^{u}(e^{u}-1)}{(e^{u}+1)^{3}u}du.
\end{align}}
{\bf Proof.} For the part sum of the series $\sum_{k=1}^{n}\frac{1}{(2k-1)^{3}}$, the following expression is valid
\begin{align}
\sum_{k=1}^{n}\frac{1}{(2k-1)^{3}}=\frac{\pi^{2}}{4}\int_{\gamma_{n}}\frac{e^{z}(e^{z}-1)}{(e^{z}+1)^{3}z}{\rm d}z,~n\in\mathbb{N}^{+},
\end{align}
where $\gamma_{n}=\Gamma_{1}^{(n)}+\Gamma_{2}^{(n)}+\Gamma_{3}^{(n)}+\Gamma_{4}^{(n)}+\Gamma_{\varepsilon}$ represent the contour with vertices $\pm R_{n}$, $ R_{n}+{\rm i}R_{n}$, $-R_{n}+{\rm i}R_{n}$ (see Fig.~1).   $\Gamma_{\varepsilon}$ is a semicircle  of radius  $\varepsilon\in(0,\pi/2)$ centered at the origin. $\Gamma_{1}^{(n)}$ represent $(-R_{n},0)\rightarrow (-\varepsilon,0)$ and $(\varepsilon,0)\rightarrow (R_{n},0)$.   $R_{n}=2n\pi$. In fact,
since in the neighborhood of the point $z=0$ and $z=z_{k}=(2k-1)\pi{\rm i}$,~$k=1,2,\cdots,n$, the following results are valid:
\begin{align}
\frac{e^{z}-1}{z}=O(1),\quad \frac{e^{z}}{(e^{z}+1)^{3}}=O((z-z_{k})^{-3}).\nonumber
\end{align}
\begin{figure}
  \centering
  \includegraphics[width=3.4in]{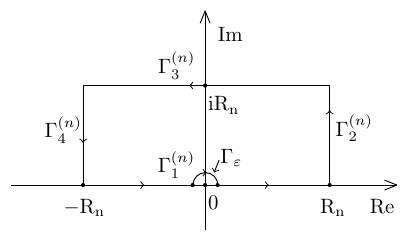}\\
 \renewcommand{\figurename}{\footnotesize{Figure}}
\caption{\footnotesize{The contour $\gamma_{n}$}}\label{BHtau20}
\end{figure}
Let $f(z)=\frac{e^{z}(e^{z}-1)}{(e^{z}+1)^{3}z}.$ The function $f(z)$ is bounded on the boundary of the square with poles at the points $z_{k}$.
It follows that the integral in $(2)$ is equal to the sum of residues of the integrand at the points $z_{1}=\pi{\rm i},z_{2}=3\pi {\rm i},\cdots,z_{n}=(2n-1)\pi{\rm i}$,
\begin{align}
\int_{\gamma_{n}}f(z){\rm d}z=2\pi {\rm i}\sum_{k=1}^{n}\underset{z=z_{k}}{\text{Res}} f(z)
=\frac{4}{\pi^{2}}\sum_{k=1}^{n}\frac{1}{(2k-1)^{3}}.
\end{align}
Furthermore, $|f(z)|=O(\frac{1}{R_{n}e^{R_{n}}})$,~$z\in \Gamma_{j}^{(n)},j=2,3,4$.
In view of $\sum\limits_{k=1}^{\infty}\frac{1}{(2k-1)^{3}}=7\zeta(3)/8$, and let $R_{n}\rightarrow +\infty$, $\varepsilon\rightarrow0$ in $\gamma_{n}$, we obtain the required expression
\begin{align}
\int_{-\infty}^{+\infty}\frac{e^{u}(e^{u}-1)}{(e^{u}+1)^{3}u}{\rm d}u=\frac{7}{2\pi^{2}}\zeta(3),
\end{align}
which leads to the desired result.
$\hfill\square$

{\bf Corollary~2.1} Remembering inverse hyperbolic secant function:
${\rm asech}(u)=\log(1/u+\sqrt{1/u^{2}-1})$,~$u\in(0,1]$, then we obtain
\begin{align}
\zeta(3)=-\frac{\pi^{2}}{14}\int_{0}^{\infty}\frac{1}{u}{\rm d} {\rm sech^{2}(u)}=\frac{\pi^{2}}{7}\int_{0}^{1}\frac{u}{{\rm asech(u)}}{\rm d}u.\nonumber
\end{align}
By using the integral transformation~$\log u=x$ and remembering that  $\int_{0}^{\infty}2(1+x)^{-3}x^{t}{\rm d}x=\pi t(1-t)/\sin(\pi t)$, $t\in (0,1)$, we find $(4)$ equals with $(22)$ in [1]:
\begin{align}
\zeta(3)=\frac{2\pi^{2}}{7}\int_{0}^{\infty}\frac{x-1}{(1+x)^{3}\log x}{\rm d}x=\frac{2\pi^{2}}{7}\int_{0}^{1}\int_{0}^{\infty}\frac{x^{t}}{(1+x)^{3}}{\rm d}x{\rm d}t=\frac{1}{7}\int_{0}^{\pi}\frac{x(\pi-x)}{\sin x}{\rm d}x,\nonumber
\end{align}
which means the result $(4)$ is a special case of the more general result in [2]. But for $\zeta(2n+1),n\geq 2$, we  found no  association between these representations.

However, let us expand the integral representation of $\zeta(3)$  into the sum of partial fractions with respect to $\exp(u)$:
\begin{align}
\zeta(3)=\frac{2\pi^{2}}{7}\sum_{l=1}^{3}\int_{-\infty}^{+\infty}\frac{w_{l}}{u(e^{u}+1)^{l}}{\rm d}u,
\end{align}
where $w_{1}=1,w_{2}=-3,w_{3}=2$.
Since $\sum\limits_{l=1}^{3}w_{l}=0$, the integrand $\sum\limits_{l=1}^{3}\frac{w_{l}}{u(e^{u}+1)^{l}}$ is an even function in $\mathbb{R}$ except a removable singularity at origin. We wonder whether the  decomposition (5) is unique and can be generalized to the general $\zeta(2n+1),n\geq2$.
\section*{3. Evaluation of $\zeta(m)$ at positive odd integers}

\subsection*{3.1 A formula for $\zeta(m),m\geq2$}

The evaluation of $\zeta(3)$ can be straightforwardly generalized to $\zeta(m)$, $m\geq2$. Before considering the general case, we need the following useful lemma.

{\bf Lemma~3.1} {\it Let $m\geq 2$ be a positive integer, and $\gamma_{n}$ represent the integral contour in Lemma 2.1.  There exists unique $w_{l}\in\mathbb{Q}$, which satisfying $\sum\limits_{l=1}^{m}w_{l}=0$, such that
\begin{align}
\sum_{l=1}^{m}\int_{\gamma_{n}}\frac{w_{l}}{z(e^{z}+1)^{l}}{\rm d}z=\frac{2S_{m}}{({\rm i}\pi)^{m-1}}\left(\zeta(m)\frac{2^{m}-1}{2^{m}}-\sum_{k=n+1}^{\infty}\frac{1}{(2k-1)^{m}}\right),~ n=1,2,\cdots,\nonumber
\end{align}
wherein $S_{m}\in \mathbb{Q}$.}

{\bf Proof.} Step {\bf 1}. Let $f_{m}(z)=\sum\limits_{l=1}^{m}\frac{w_{l}}{z(e^{z}+1)^{l}}$. Note that the integrand $f_{m}(z)$ over the vertical and top side  $\Gamma^{(n)}_{j}$, $j=2,3,4$ tend to $0$, when $R_{n}\rightarrow\infty$, since it can be bounded by $Ce^{-R_{n}}/R_{n}$ and constant $C>0$. Let $z_{k}=(2k-1)\pi{\rm i}$,~$k=1,2,\cdots,n$. By means of the theorem of residues, we have
\begin{align}
\int_{\gamma_{n}}f_{m}(z){\rm d}z\nonumber
&=2\pi {\rm i}\sum_{k=1}^{n}\lim_{z\rightarrow z_{k}}\frac{1}{(m-1)!}\left(\frac{{\rm d}}{{\rm d}z}\right)^{m-1}\sum_{l=1}^{m}\frac{w_{l}(z-z_{k})^{m}}{z(e^{z}+1)^{l}}\\\nonumber
&=2\pi {\rm i}\sum_{k=1}^{n}\lim_{z\rightarrow z_{k}}\frac{1}{(m-1)!}\sum_{l=1}^{m}w_{l}\sum_{j=0}^{m-1}\binom{m-1}{j}\left[\frac{(z-z_{k})^{l}}{z(e^{z}+1)^{l}}\right]^{(j)}\left[(z-z_{k})^{m-l}\right]^{(m-j-1)}\\
&=2\pi {\rm i}\sum_{k=1}^{n}\sum_{l=1}^{m}\frac{w_{l}}{(l-1)!}\lim_{z\rightarrow z_{k}}\left[\frac{(z-z_{k})^{l}}{z(e^{z}+1)^{l}}\right]^{(l-1)}.
\end{align}

Step {\bf 2}. Let us denote
\begin{align}
g_{l}(z)=\frac{(z-z_{k})^{l}}{z(e^{z}+1)^{l}},~
D_{l}=\lim_{z\rightarrow z_{k}}\frac{g_{l}^{(l-1)}(z)}{(l-1)!},~l=1,2,\cdots,m.
\end{align}
The conclusions
\begin{equation*}
g_{l}^{(l-1)}(z)=\sum_{j=0}^{l-1}\binom{l-1}{j}\left(\frac{1}{z}\right)^{(l-1-j)}\left[\left(\frac{z-z_{k}}{e^{z}+1}\right)^{l}\right]^{(j)},
\end{equation*}
follow from the Leibniz rule for the differentiation of a product. For $j=0$ in $g_{l}^{(l-1)}(z)$, we have
\begin{align}
\lim_{z\rightarrow z_{k}}\binom{l-1}{0}\left(\frac{1}{z}\right)^{(l-1)}\left(\frac{z-z_{k}}{e^{z}+1}\right)^{l}=-\frac{(l-1)!}{z_{k}^{l}}.
\end{align}
According with $\big[\frac{1}{(e^{z}+1)^{l}}\big]'+\frac{l}{(e^{z}+1)^{l}}=\frac{l}{(e^{z}+1)^{l+1}}$, which means for all $l$, the residues $c_{-1}^{(l)}$ of $(e^{z}+1)^{-l}$ at $z=z_{k}$ are equal:
\begin{equation*}
c_{-1}^{(l)}=c_{-1}^{(1)}=\frac{1}{2\pi {\rm i}}\int_{\Lambda_{k}}\frac{1}{(e^{z}+1)^{l}}{\rm d}z,~l=1,2,\cdots,m,
\end{equation*}
where $\Lambda_{k}$ is the small circle centered at the poles $z_{k}$ and of radius $\epsilon_{k}>0$. When $l=1$, we acquire $c_{-1}^{(l)}=c_{-1}^{(1)}=\lim\limits_{z\rightarrow z_{k}}(z-z_{k})/(e^{z}+1)=-1$. For $j=l-1$ in $g_{l}^{(l-1)}(z)$, by using Laurent  series of  $(e^{z}+1)^{-l}$ at $z=z_{k}$, we obtain
\begin{align}
\lim_{z\rightarrow z_{k}}\binom{l-1}{l-1}\frac{1}{z}\left[\left(\frac{z-z_{k}}{e^{z}+1}\right)^{l}\right]^{(l-1)}=\lim_{z\rightarrow z_{k}}\frac{1}{z}[c_{-1}^{(l)}(z-z_{k})^{l-1}]^{(l-1)}=-\frac{(l-1)!}{z_{k}}.
\end{align}
As a result of the relations $(8)$ and $(9)$, we can  write out $D_{l}$:
\begin{align}
D_{l}=-(\frac{b_{l}^{(l)}}{z_{k}^{l}}+\cdots+\frac{b_{j}^{(l)}}{z_{k}^{j}}+\cdots+\frac{b_{1}^{(l)}}{z_{k}}),~j=1,\cdots,l,
\end{align}
in which $b_{l}^{(l)}=b_{1}^{(l)}=1$, $b_{2}^{(l)},\cdots,b_{l-1}^{(l)}\in \mathbb{Q}$.

Step {\bf 3}. In this part, we will prove that there exists unique $w_{l}$, which satisfying
\begin{equation*}
\sum_{l=1}^{m}w_{l}D_{l}=\frac{S_{m}}{z_{k}^{m}},~S_{m}\in \mathbb{Q}.
\end{equation*}
In fact, it should be noted that, by applying $(10)$ in Step {\bf 2}, we obtain
\begin{align}
\sum_{l=1}^{m}w_{l}D_{l}=-\sum_{l=1}^{m}w_{l}(\frac{1}{z_{k}^{l}}+\cdots+\frac{b_{j}^{(l)}}{z_{k}^{j}}+\cdots+\frac{1}{z_{k}})=\frac{S_{m}}{z_{k}^{m}}.
\end{align}
Comparing the coefficients on both sides of identity $(11)$ with respect to $1/z^{l}_{k}$,  we have the relation:
\begin{align}
\begin{pmatrix}
1&1&\cdots&1\\
0&1&\cdots&b_{2}^{(m)}  \\
\vdots&\vdots&\ddots&\vdots \\
0&0&\cdots&1
\end{pmatrix}_{m\times m}
\begin{pmatrix}
w_{1}\\
w_{2}\\
\vdots\\
w_{m}
\end{pmatrix}=
\begin{pmatrix}
0\\
0\\
\vdots\\
-S_{m}
\end{pmatrix}
.
\end{align}
Let $A_{m\times m}$ represents the coefficient matrix of $(12)$. Since $A_{m\times m}$ is an upper-triangular matrix, and $\|A_{m\times m}\|=1,$ then by Cramer's Rule, the existence and unique of $w_{l}$ can be obtained. Set $C^{(m)}=[0,0,\cdots,-S_{m}]^{T}$, $W^{(m)}=[w_{1},w_{2},\cdots,w_{m}]^{T}$, then
\begin{align} w_{l}=[W^{(m)}]_{l},~l=1,\cdots,m,
\end{align}
in which $W^{(m)}=A_{m\times m}^{-1}C^{(m)}$ and $[\cdot]_{l}$ represents the $l-$th element in the column vector.

Given that $z_{k}=(2k-1)\pi{\rm i}$, and substitute $w_{l}$ and $D_{l}$ into $(6)$,  we obtain
\begin{equation*}
\begin{aligned}
\sum_{l=1}^{m}\int_{\gamma_{n}}\frac{w_{l}}{z(e^{z}+1)^{l}}{\rm d}z=2\pi {\rm i}\sum_{k=1}^{n}\frac{S_{m}}{z_{k}^{m}}=\frac{2S_{m}}{({\rm i}\pi)^{m-1}}\Big[\zeta(m)(1-\frac{1}{2^{m}})-\sum_{k=n+1}^{\infty}\frac{1}{(2k-1)^{m}}\Big],
\end{aligned}
\end{equation*}
which gives the required equivalence. Lemma 3.1 is thus proved.
$\hfill\square$

{\bf Corollary~3.1}
Let
\begin{equation*}
S_{m}=\left\{
\begin{aligned}
&(-1)^{\frac{m-1}{2}}\Gamma(m),~~m=3,5,\cdots,\\
&(-1)^{\frac{m+2}{2}}\Gamma(m),~~m=2,4,\cdots,
\end{aligned}
\right.
\end{equation*}
where $\Gamma(\cdot)$ denotes the Gamma function and $\Gamma(m)=(m-1)!$. For $n=1,2,\cdots$, by using Lemma 3.1, there exists unique $w_{l}\in \mathbb{Q}$, $l=1,2,\cdots,m$, such that
\begin{align}
\sum_{l=1}^{m}\int_{\gamma_{n}}\frac{w_{l}}{z(e^{z}+1)^{l}}{\rm d}z=
\frac{2\Gamma(m)}{\pi^{m-1}{\rm i}^{\kappa}}\Big[\zeta(m)(1-\frac{1}{2^{m}})-\sum_{k=n+1}^{\infty}\frac{1}{(2k-1)^{m}}\Big],
\end{align}
in which $\kappa=0$ or $1$ with respect to $m$ is positive odd or even integer.

\subsection*{3.2 A formula for $\zeta(m),m\geq3$ at positive odd integers  }

The generalized Bernoulli polynomials[4] $B_{n}^{(l)}(x)$ of degree $n$ in $x$ are defined by the generating function:
\begin{align}
\left(\frac{z}{e^{z}-1}\right)^{l}e^{xz}=\sum_{n=0}^{\infty}B_{n}^{(l)}(x)\frac{z^n}{n!},~|z|<2\pi,1^{l}:=1,
\end{align}
for arbitrary (real or complex) parameter $l$ and $B_{n}^{(l)}(0)=B_{n}^{(l)}$.
Let $m\geq 3$ be a positive odd integer. $w_{l}$ are solutions to $(12)$ and $z_{k}=(2k-1)\pi{\rm i}$, $k\geq1$. Set $x=z-z_{k}$ and $0<|x|<2\pi$. Using the property (22) of $B_{n}^{(l)}(x)$ in [4],  we obtain
\begin{align}
\frac{1}{(e^{z}+1)^{l}}&=\frac{(-1)^{l}}{(e^{x}-1)^{l}}=\sum_{n=0}^{\infty}(-1)^{l}B_{n}^{(l)}\frac{x^{n-l}}{n!},\\
\frac{1}{(e^{-z}+1)^{l}}&=\frac{e^{xl}}{(e^{x}-1)^{l}}=\sum_{n=0}^{\infty}B_{n}^{(l)}(l)\frac{x^{n-l}}{n!}=\sum_{n=0}^{\infty}(-1)^{n}B_{n}^{(l)}\frac{x^{n-l}}{n!}.
\end{align}
Since the Laurent  series of  $(e^{z}+1)^{-l}$ centered at $z=z_{k}$ have the form
\begin{align}
\frac{1}{(e^{z}+1)^{l}}=\sum_{j=-l}^{\infty}c_{j}^{(l)}(z-z_{k})^{j},~l=1,2,\cdots,m,
\end{align}
where
\begin{align}
c_{j}^{(l)}=\frac{1}{2\pi i}\int_{\Lambda_{k}}\frac{{\rm d}z}{(e^{z}+1)^{l}(z-z_{k})^{j+1}},~j=-l,-(l-1),\cdots.\nonumber
\end{align}
Comparing the coefficients of $1/(z-z_{k})^{l}$ in $(16)$ and $(18)$, we obtain
\begin{align}
c_{j}^{(l)}=\frac{(-1)^{l}B_{l+j}^{(l)}}{\Gamma(l+j+1)},~j=-l,\cdots,-1.
\end{align}
Taking $(7)$ and $(18)$ into account, we obtain
\begin{align}
D_{l}=\sum_{j=0}^{l-1}(-1)^{l-j-1}\binom{l-1}{j}\frac{\Gamma(l-j)\Gamma(j+1)c_{-(l-j)}^{(l)}}{\Gamma(l)z_{k}^{l-j}}
=\sum_{j=1}^{l}(-1)^{j+1}\frac{c_{-j}^{(l)}}{z_{k}^{j}}.\nonumber
\end{align}
Hence, it follows from $(10)$ that
\begin{align}
b_{j}^{(l)}=(-1)^{j}c_{-j}^{(l)}.
\end{align}
Combining with $(19)$, we deduce that
\begin{align}
b_{j}^{(l)}=\frac{(-1)^{l-j}}{\Gamma(l-j+1)}B_{l-j}^{(l)},~j=1,2,\cdots,l.
\end{align}
With the help of $(17)$ and $(12)$, we obtain
\begin{align}
\begin{aligned}
\frac{1}{2\pi {\rm i}}\sum_{l=1}^{m}\int_{\gamma_{n}}\frac{w_{l}}{z(e^{-z}+1)^{l}}{\rm d}z&=\frac{1}{2\pi{\rm i}}\sum_{k=1}^{n}\sum_{j=1}^{m}\sum_{l=j}^{m}\int_{\Lambda_{k}}w_{l}(-1)^{l-j}\frac{B_{l-j}^{(l)}}{\Gamma(l-j+1)}\frac{1}{z(z-z_{k})^{j}}{\rm d}z\\
&=\sum_{k=1}^{n}\sum_{j=1}^{m}\sum_{l=j}^{m}w_{l}b_{j}^{(l)}\frac{(-1)^{j-1}}{z_{k}^{j}}=\sum_{k=1}^{n}\frac{(-1)^{m-1}w_{m}b_{m}^{(m)}}{z_{k}^{m}}\\
&=-\sum_{k=1}^{n}\frac{S_{m}}{z_{k}^{m}}.
\end{aligned}\nonumber
\end{align}
Hence, we have the relation
\begin{align}
\sum_{l=1}^{m}\int_{\gamma_{n}}\frac{w_{l}}{z(e^{z}+1)^{l}}{\rm d}z&=\frac{1}{2}\sum_{l=1}^{m}\int_{\gamma_{n}}\frac{w_{l}}{z}\left[\frac{1}{(e^{z}+1)^{l}}-\frac{1}{(e^{-z}+1)^{l}}\right]{\rm d}z\nonumber\\
&=\frac{1}{2}\sum_{l=1}^{m}\sum_{p=0}^{l-1}\int_{\gamma_{n}}\frac{(1-e^{z})w_{l}e^{zp}}{z(e^{z}+1)^{l}}{\rm d}z.
\end{align}
Since $z=0$ is a removable singularity pole for the integrand in $(22)$, let $R_{n}\rightarrow +\infty$ and $\varepsilon\rightarrow 0$ in $\gamma_{n}$, then we have
\begin{align}
2\Gamma(m)(1-\frac{1}{2^{m}})\frac{\zeta(m)}{\pi^{m-1}}=\sum_{l=1}^{m}\int_{-\infty}^{+\infty}\frac{w_{l}}{u(e^{u}+1)^{l}}{\rm d}u=\frac{1}{2}\sum_{l=1}^{m}\sum_{p=0}^{l-1}\int_{-\infty}^{+\infty}\frac{w_{l}(1-e^{u})e^{up}}{u(e^{u}+1)^{l}}{\rm d}u
.\nonumber
\end{align}
Note that the integrand $\sum\limits_{l=1}^{m}\sum\limits_{p=0}^{l-1}\frac{w_{l}(1-e^{u})e^{up}}{u(e^{u}+1)^{l}}$ is an even function on $\mathbb{R}$ except a removable singularity at origin, then we have
\begin{align}
\zeta(m)=\frac{(2\pi)^{m-1}}{(2^{m}-1)\Gamma(m)}\sum_{l=1}^{m}\sum_{p=0}^{l-1}\int_{0}^{+\infty}\frac{w_{l}(1-e^{u})e^{up}}{u(e^{u}+1)^{l}}{\rm d}u, ~~m=3,5,\cdots,
\end{align}
in which $w_{l}$ are determined by $(13)$.

\section*{4. Evaluation of $\zeta(5)$ and $\zeta(7)$ by using the generalized Bernoulli numbers }

In this part, we calculate the numerical results of $\zeta(5)$ and $\zeta(7)$. In view of the Srivastava-Todorov formula [8, p.510, Eq.(3)], we know the following representation for the generalized Bernoulli numbers:
\begin{align}
B_{n}^{(l)}=\sum_{k=0}^{n}\binom{l+n}{n-k}\binom{l+k-1}{k}\frac{n!}{(n+k)!}\sum_{j=0}^{k}(-1)^{j}\binom{k}{j}j^{n+k}.
\end{align}
Using the formula $(21)$ and $(24)$, it is verified by a straightforward calculation that
\begin{equation*}
\begin{aligned}
&D_{1}=-\frac{1}{z_{k}},
~D_{2}=-\frac{1}{z_{k}^{2}}-\frac{1}{z_{k}},
~D_{3}=-\frac{1}{2!}(\frac{2}{z_{k}^{3}}+\frac{3}{z_{k}^{2}}+\frac{2}{z_{k}}),\\
&D_{4}=-\frac{1}{3!}(\frac{3!}{z_{k}^{4}}+\frac{12}{z_{k}^{3}}+\frac{11}{z^{2}_{k}}+\frac{3!}{z_{k}}),
~D_{5}=-\frac{1}{4!}(\frac{4!}{z^{5}_{k}}+\frac{60}{z_{k}^{4}}+\frac{70}{z_{k}^{3}}+\frac{50}{z_{k}^{2}}+\frac{4!}{z_{k}}),\\
&D_{6}=-\frac{1}{5!}(\frac{5!}{z_{k}^{6}}+\frac{360}{z_{k}^5}+\frac{510}{z_{k}^4}+\frac{450}{z_{k}^3}+\frac{274}{z_{k}^2}+\frac{5!}{z_{k}}),\\
&D_{7}=-\frac{1}{6!}(\frac{6!}{z_{k}^{7}}+\frac{2520}{z_{k}^{6}}+\frac{4200}{z_{k}^{5}}+\frac{4410}{z_{k}^{4}}+\frac{3248}{z_{k}^{3}}+\frac{1764}{z_{k}^2}+\frac{6!}{z_{k}}),
\end{aligned}
\end{equation*}
where $z_{k}=(2k-1)\pi{\rm i}, 1\leq k\leq n$.
According with $S_{m}$ in Corollary $3.1$ and solving linear equation $(12)$, we obtain $w_{l}$ and the following result.\\
{\bf Theorem 4.1}{\it
~~~{\rm (1):} For Riemann zeta function value $\zeta(5)$
\begin{align}
\zeta(5)=\frac{(2\pi)^{4}}{(2^{5}-1)\Gamma(5)}\sum_{l=1}^{5}\int_{-\infty}^{\infty}\frac{w_{l}}{(e^{u}+1)^{l}u}{\rm d}z,\nonumber
\end{align}
where $w_{1}=-1,w_{2}=15,w_{3}=-50,w_{4}=60,w_{5}=-24$.

{\rm (2):} For Riemann zeta function value $\zeta(7)$
\begin{align}
\zeta(7)=\frac{(2\pi)^{6}}{(2^{7}-1)\Gamma(7)}\sum_{l=1}^{7}\int_{-\infty}^{\infty}\frac{w_{l}}{(e^{u}+1)^{l}u}{\rm d}z,\nonumber
\end{align}
where $w_{1}=1,w_{2}=-63,w_{3}=602,w_{4}=-2100,w_{5}=3360,w_{6}=-2520,w_{7}=720$.
}

\section*{5. The linear combination of $\zeta(2n+1)/\pi^{2n},n\geq1$ }

In this part, we investigate the integral representation for the linear combination of $\zeta(2n+1)/\pi^{2n}$. For all positive odd integer $m\geq 3$, it follows from $(23)$ that
\begin{align}
2\Gamma(m)(1-\frac{1}{2^{m}})\frac{\zeta(m)}{\pi^{m-1}}=-\sum_{l=1}^{m}\sum_{p=0}^{l-1}\int_{0}^{+\infty}\frac{w_{l}e^{(2p+1-l)u}}{u(e^{u}+e^{-u})^{l}}{\rm d}(e^{u}+e^{-u}).
\end{align}

Considering the symmetry of the power of $\exp(2p+1-l)$ with $p$ from $0$ to $l-1$, we can draw the following conclusions:
\begin{align}
\sum_{p=0}^{l-1}\frac{e^{(2p+1-l)u}}{(e^u+e^{-u})^{l}}=
\sum_{j=1}^{k}\frac{q_{j}^{(l)}}{(e^{u}+e^{-u})^{2j-1}},~~k=\lceil{l/2}\rceil,q_{j}^{(l)}\in\mathbb{Z},
\end{align}
in which $\lceil{\cdot}\rceil$ means the ceiling function and $q_{j}^{(l)}=1-\sum\limits_{\kappa=1}^{j-1}\binom{l+1-2\kappa}{j-\kappa}q_{\kappa}^{(l)},$ for $j\geq 2$, $q_{1}^{(l)}\equiv1$.
Remembering that $\sum\limits_{l=1}^{m}w_{l}=0$ in Lemma 3.1, then we obtain
\begin{align}
2\Gamma(m)(1-\frac{1}{2^{m}})\frac{\zeta(m)}{\pi^{m-1}}&=-\int_{0}^{\infty}\sum_{j=1}^{\lceil{m/2}\rceil}[\sum_{l=2j-1}^{m}w_{l}q_{j}^{(l)}]\frac{{\rm d}(e^{u}+e^{-u})}{u(e^{u}+e^{-u})^{2j-1}}\nonumber\\
&=\int_{0}^{\infty}\sum_{j=2}^{\lceil{m/2}\rceil}[\sum_{l=2j-1}^{m}w_{l}q_{j}^{(l)}]\frac{2^{-2j+1}}{(j-1)u}{\rm d}({\rm sech}{(u)})^{2j-2}\nonumber\\
&=-\int_{0}^{1}\sum_{j=2}^{\lceil{m/2}\rceil}\frac{1}{4^{j-1}}[\sum_{l=2j-1}^{m}w_{l}q_{j}^{(l)}]\frac{u^{2j-3}}{{\rm~asech}{(u)}}{\rm d}u.
\end{align}
Let us denote
\begin{align}
\tau_{j}^{(m)}=-\frac{2^{m-1}}{\Gamma(m)(2^{m}-1)}\sum\limits_{l=2j-1}^{m}\frac{w_{l}q_{j}^{(l)}}{4^{j-1}},
\end{align}
then
\begin{align}
\frac{\zeta(m)}{\pi^{m-1}}=\sum_{j=2}^{\lceil{m/2}\rceil}\tau_{j}^{(m)}\int_{0}^{1}\frac{u^{2j-3}}{{\rm~asech}{(u)}}{\rm d}u.
\end{align}
Through above discussion, we have the following theorem.

{\bf Theorem 5.1} {\it For any positive integer $n$,
then there exists $\theta_{k}\in\mathbb{Q}$, $k=1,2\cdots,n$, so that
\begin{equation*}
\sum_{k=1}^{n}\theta_{k}\frac{\zeta(2k+1)}{\pi^{2k}}=
\theta_{n+1}I_{n},~\theta_{n+1}=0 ~or ~1,
\end{equation*}
where $I_{n}$ satisfy
\begin{equation*}
I_{n}=\int_{0}^{1}\frac{u^{2n-1}}{{\rm asech}(u)}du>0,~~\lim_{n\rightarrow\infty}I_{n}=0.
\end{equation*}}
{\bf Proof.} Since ${\rm asech(u)}=\log(1/u+\sqrt{1/u^{2}-1})$, then for $\vartheta$ near to $0$ from right side
\begin{align}
\frac{1}{{\rm asech(u)}}=-\frac{1}{\log u}-\frac{\log2}{\log^{2}u}+o(1/\log^{2}u),~0<u<\vartheta<1,
\end{align}
where $o(\cdot)$ means higher-order infinitesimal. Hence, we obtain $1/{\rm asech}(u)\rightarrow0$ as $u\rightarrow 0^{+}$.
And for $\eta>0$ near to $1$ from left side, when  $u\in(\eta,1)$,
\begin{align}
\frac{1}{{\rm asech(u)}}=\frac{\sqrt{2}}{2\sqrt{1/u-1}}+\frac{\sqrt{2}}{24}\sqrt{1/u-1}+o(\sqrt{1/u-1}),
\end{align}
then ${\rm asech(u)}/\sqrt{1/u-1}\rightarrow\sqrt{2}$ as $u\rightarrow 1^{-}$. Since $\int_{0}^{1}(1/u-1)^{-1/2}du=\pi/2$,
so we can deduce that $\int_{0}^{1}1/{\rm asech(u)}{\rm du}$ converges.  According with ${\rm asech(u)}>0$, when $u\in(0,1)$, and $0<I_{n}\leq\int_{0}^{1}1/{\rm asech(u)}{\rm du}$, then $\{I_{n}\}_{n\geq1}$ is a monotonically decreasing positive sequence with a lower bound, so it must converge.

From $(31)$, when $u\in(\eta,1)$, we can find proper $M>0$, so that
${\rm 1/asech(u)}\leq \frac{M}{\sqrt{1/u-1}}$. As a result, we have
\begin{align}
I_{n}&=\int_{0}^{1}\frac{u^{2n-1}}{{\rm asech}(u)}du,\nonumber\\
&=\int_{0}^{\eta}\frac{u^{2n-1}}{{\rm asech}(u)}du+\int_{\eta}^{1}\frac{u^{2n-1}}{{\rm asech}(u)}du\nonumber\\
&\leq\eta^{2n-1}\int_{0}^{1}\frac{1}{{\rm asech}(u)}du+M\int_{0}^{1}u^{2n-1}(1/u-1)^{-1/2}du\nonumber\\
&\leq\eta^{2n-1}\int_{0}^{1}\frac{1}{{\rm asech}(u)}du+M\sqrt{\pi}\frac{\Gamma(2n+1/2)}{\Gamma(2n+1)},\nonumber
\end{align}
since $\eta\in(0,1)$ and using Stirling's approximation: $\frac{\Gamma(2n+1/2)}{\Gamma(2n+1)}=O(1/\sqrt{2n+1})$ as $n\rightarrow\infty$, so $\lim\limits_{n\rightarrow\infty}I_{n}=0$.

Using $(29)$, let $m=2n+1$, if for $\forall n\geq1$, then $\tau_{n+1}^{(2n+1)}\neq0$, we can choose proper $\theta_{k}$, $1\leq k\leq n$ so that $\theta_{n+1}= 1$. If there exists $n\geq2$ such that $\tau_{n+1}^{(2n+1)}=0$, then we can choose proper $\theta_{k}$ (where $\theta_{k}$ satisfying $\sum_{k=1}^{n}\theta_{k}^{2}\neq0$) making $\theta_{n+1}=0$.
Therefore, there must exist proper $\theta_{k}\in \mathbb{Q}$ which make $\sum\limits_{k=1}^{n}\theta_{k}\frac{\zeta(2k+1)}{\pi^{2k}}=\theta_{n+1}I_{n}$, $\theta_{n+1}=0$ or $1$. Thus Theorem 5.1 is proved.  $\hfill\square$

Since the problem of the irrationality of the values of the Riemann zeta function at odd integers is one of the most attractive topics. Before Apery's results[10], the arithmetic properties of the zeta function at odd points seemed inaccessible.
According with Theorem $5.1$ in this paper and Nesterenko's theorem[9], we derive criteria for the dimension of the vector space spanned over the rational by $\zeta(2n+1)/\pi^{2n}$, $n\geq1$.

{\bf Theorem 5.2} {\it
Let $m=2n+1$ in $(29)$, if  $\tau_{n+1}^{(2n+1)}\neq0$ in $(28)$ for $n\geq1$, then
\begin{align}
\dim_{\mathbb{Q}}(\mathbb{Q}\frac{\zeta(3)}{\pi^{2}}+\mathbb{Q}\frac{\zeta(5)}{\pi^{4}}+\cdots+\mathbb{Q}\frac{\zeta(2n+1)}{\pi^{2n}}+\cdots)=\infty;
\end{align}
If there exits $n\geq2$, so that $\tau_{n+1}^{(2n+1)}=0$, then
\begin{align}
\dim_{\mathbb{Q}}(\mathbb{Q}\frac{\zeta(3)}{\pi^{2}}+\mathbb{Q}\frac{\zeta(5)}{\pi^{4}}+\cdots+\mathbb{Q}\frac{\zeta(2n+1)}{\pi^{2n}})\leq n-1.
\end{align}
}
\section*{6. Concluding remarks }
Through numerical  calculations  of $\zeta(5)$, $\zeta(7)$ and $\zeta(9)$, we do believe that only conclusion (32) is possible, which means the sequence $\zeta(3)/\pi^{2},\zeta(5)/\pi^{4},\zeta(7)/\pi^{6}$,$\cdots$ contains infinitely many irrational numbers.
But it's a pity that we can't prove it theoretically. This is also one of our future research contents.

\section*{Conflict of interest}

The author declares that there is no conflict of interest with respect to the publication of this paper.

\section*{Funding}

Work funded by High-level Talent Research Project Fund of NHU (Grant Nos. 2022rcjj24).

\end{document}